\documentclass[12pt]{amsart}
\usepackage{latexsym, amsmath,amssymb}
\usepackage{amsthm}
\usepackage{hyperref}
\usepackage{xcolor}
\hypersetup{
    colorlinks,
    linkcolor={red!50!black},
    citecolor={blue!50!black},
    urlcolor={blue!80!black}
}

\usepackage{graphicx}

\newtheorem*{theorem*}{Theorem}

 \numberwithin{equation}{section}

\def\XXint#1#2#3{{\setbox0=\hbox{$#1{#2#3}{%
\int}$ }
\vcenter{\hbox{$#2#3$ }}\kern-.6\wd0}}

\renewcommand{\epsilon}{\varepsilon}
\newtheorem{theorem}{Theorem}

\newtheorem{lemma}[theorem]{Lemma}
\newtheorem{corr}[theorem]{Corollary}

\newtheorem{proposition}[theorem]{Proposition}
\newtheorem{deff}[theorem]{Definition}

\newtheorem{conjecture}[theorem]{Conjecture}
\newcommand{\bth}{\begin{theorem}}
\newcommand{\ble}{\begin{lemma}}
\newcommand{\bcor}{\begin{corr}}

\newcommand{\bdeff}{\begin{deff}}

\newcommand{\bprop}{\begin{proposition}}
\newcommand{\ele}{\end{lemma}}
\newcommand{\ecor}{\end{corr}}
\newcommand{\edeff}{\end{deff}}

\numberwithin{theorem}{section}

\newcommand{\eprop}{\end{proposition}}

\newcommand{\Rn}{{\mathbb R}^n}

\renewcommand{\Pi}{\varPi}

\renewcommand{\epsilon}{\varepsilon}

\begin{document}

\title[Distribution symmetry of toral eigenfunctions]
{Distribution symmetry\\ of toral eigenfunctions}

\author[A. D. Mart\'inez]{\'Angel D. Mart\'inez} \address{Albacete, Fields Ontario Postdoctoral Fellow, University of Toronto Mississauga, Canad\'a} \email{amartine@fields.toronto.edu}
\author[F. Torres de Lizaur]{Francisco Torres de Lizaur} \address{Sevilla, Fields Ontario Postdoctoral Fellow, University of Toronto, Canad\'a}  \email{fj.torres@icmat.es}

\maketitle

\begin{abstract}
In this paper we study a number of conjectures on the behavior of the value distribution of eigenfunctions. On the two dimensional torus we observe that the symmetry conjecture holds in the strongest possible sense. On the other hand we provide a counterexample for higher dimensional tori, which relies on a computer assisted argument. Moreover we prove a theorem on the distribution symmetry of a certain class of trigonometric polynomials that might be of independent interest. 
\end{abstract}

\section{Introduction}

The importance of Laplace eigenfunctions might be underpinned by the fact that they are analogous to the trigonometrical polynomials in the classical harmonic analysis. How far the analogy can be drawn lies at the heart of many conjectures on their behavior. For instance, Yau conjectured that the nodal set of Laplace eigenfunctions, $-\Delta_g\psi_{\lambda}=\lambda\psi_{\lambda}$, on a smooth Riemannian manifold $M$ of dimension $n$, has $(n-1)$-dimensional Hausdorff measure comparable to $\lambda^{1/2}$. In the case of eigenfunctions on the torus this is easily seen to be the case as an application of the fundamental theorem of calculus and an elementary geometric argument. For real-analytic manifolds, the proof of Yau's conjectue was given by Donnelly and Fefferman in \cite{DF}. Recently, Logunov proved the lower bound in the smooth category and improved the exponential upper bound of Hardt and Simon to a polynomial bound (cf. \cite{Lo, Lo2}). It is beyond the scope of this paper to provide a comprehensive introduction to the subject and we refer the reader to the extensive literature (see e.g \cite{Z} and references therein).
 
In view of the predictions of the random wave conjectures of quantum chaos \cite{B, JNT, KRud}, it is interesting to investigate the relationship between positive and negative parts of real eigenfuntions on Riemannian manifolds. In the aforementioned seminal work, Donelly and Fefferman also proved

\begin{theorem*}[Corollary 7.10 in \cite{DF}]
There exist a constant $C>0$ such that
\[\frac{1}{C}\leq \frac{\operatorname{vol}(\{x\in M:\psi_{\lambda}(x)>0\})}{\operatorname{vol}(\{x\in M:\psi_{\lambda}(x)<0\})}\leq C.\]
\end{theorem*}
The constant depends on the manifold but not on the eigenvalue. In the case of the sphere, this quasi-symmetry result was conjectured in \cite{ArG}.  In the case of surfaces, a different proof was found by Nadirashvili in \cite{N}; local versions have appeared for instance in the work of Nazarov, Polterovich and Sodin \cite{NPS} . For smooth Riemannian manifolds the analogous quasi-symmetry statement remains open (even in two dimensions).

In connection with this, it has even been conjectured that

\begin{conjecture}[Symmetry]\label{conj1}
The limit
\[ \frac{\operatorname{vol}(\{x\in M:\psi_{\lambda}(x)>0\})}{\operatorname{vol}(\{x\in M:\psi_{\lambda}(x)<0\})}\rightarrow 1\]
holds as $\lambda$ grows to infinity.
\end{conjecture}

 An heuristic justification of the above, in the case of metrics of negative curvature, seems to be provided by Berry's conjecture, but Conjecture \ref{conj1} is actually believed to hold for any manifold \cite{LoSB}. 
 
 In this paper we will disprove Conjecture \ref{conj1} in the case of the $n$-dimensional tori with $n\geq 3$, although it will be proved to hold in $\mathbb{T}^2$ (cf. Theorem \ref{cx} and Theorem \ref{sign} below respectively). It is not clear to the authors whether other two dimensional manifolds will satisfy the conjecture. This is a particularly interesting geometry to study  the behaviour of eigenfunctions with subtle connections to number theory (cf. \cite{J, RL, RW}). In a forthcoming work the authors will prove that, for the two dimensional sphere, the conjecture holds for a basis of eigenfunctions (this is a weaker statement that trivially holds for $n$-dimensional tori regardless of the dimension).

Another well-known result explores the ratio of global extrema 

\begin{theorem*}[Nadirashvili, \cite{N}]
There exist a constant $C>0$
such that
\[\frac{1}{C}\leq \frac{\|\psi_{\lambda}\chi_{\{\psi_{\lambda}>0\}}\|_{\infty}}{\|\psi_{\lambda}\chi_{\{\psi_{\lambda}<0\}}\|_{\infty}} \leq C.\]
\end{theorem*}

In general the constant cannot  be taken to be one (which would be optimal) as shown by Jakobson and Nadirashvili using zonal spherical harmonics in \cite{JN}. This theorem was extended by Jakobson and Nadirashvili for general $L^p$ norms, $1\leq p<\infty$ (loc. cit.). The exact value of the constant for the sphere $\mathbb{S}^n$ was considered in a work of Armitage in \cite{Ar}. In \cite{JNT} Jakobson, Nadirashvili and Toth claim: \textit{it is unclear  whether }

\begin{conjecture}\label{conj2}
\textit{$\|\psi_{\lambda}\chi_{\{\psi_{\lambda}>0\}}\|_{p}/\|\psi_{\lambda}\chi_{\{\psi_{\lambda}<0\}}\|_{p} \rightarrow 1$ as $\lambda\rightarrow \infty$ for $1<p<\infty$ on a given manifold.}
\end{conjecture}

One of the by-products of the results in this paper will be to answer this affirmatively for a wide class of two dimensional tori.  It is quite likely that our counterexamples (Theorem \ref{sign}) to Conjecture \ref{conj1} will also disprove this in the case of higher dimensional tori (for instance, it is easy to observe that they do disprove the endpoint $p=\infty$), but we will not pursue that question in this paper. 

\section{Statements of results} 

The first  observation  of this paper will be the following 

 \begin{theorem}[Sign equidistribution]\label{sign}
Given a non constant real eigenfunction $\psi$ of the flat two dimensional torus, the following identity holds
\[\operatorname{vol}(\{x\in\mathbb{T}^2:\psi(x)>0\})=\operatorname{vol}(\{x\in\mathbb{T}^2:\psi(x)<0\}).\]
\end{theorem}

This is stronger than  Conjecture \ref{conj1}  and provides the first example for which it holds (to the best of the  authors'   knowledge). Let us observe in passing though, that  one can show the symmetry conjecture holds for a canonical basis of eigenfunctions for the Dirichlet problem on a ball using Stoke's aproximations of the zeroes of Bessel functions  (cf. \cite{Watson}, p. 505). This might provide an idea of the intrinsic analytic difficulties even in particular well-known cases.

The following also holds: 

\begin{theorem}[Global extrema of eigenfunctions]\label{maxmin}
Given a non constant real eigenfunction $\psi$ of the flat two dimensional torus, the following identity holds
\[\max_{x\in\mathbb{T}^2}\psi(x)=-\min_{x\in\mathbb{T}^2}\psi(x).\]
Or, stated differently, the absolute values of the maxima and minima coincide.
\end{theorem}

This fact is in clear contrast with the result of Armitage on the sphere  \cite{Ar} where equality is shown to fail. Both observations suggest some symmetry of the distribution function and will be consequences of the following more general 

\begin{theorem}[Distributional symmetry]\label{main}
Given a non constant real eigenfunction $\psi$ of the flat two dimensional torus, the distribution function
\[
d\lambda(s)=d\operatorname{vol}(\{x\in\mathbb{T}^2:\psi(x)>s\})
\]
is symmetric around $s=0$.
\end{theorem}

This might be compared with the recent work of Klartag \cite{Klartag}. 

\begin{corr}[$L^p$ norm symmetry]\label{lp}
Under the same hypothesis, the following holds
\[\int_{\{\psi>0\}}|\psi(x)|^pdx=\int_{\{\psi<0\}}|\psi(x)|^pdx\]
\end{corr}

This answers affirmatively the question raised by Jakobson, Nadirashvili and Toth (i.e. Conjecture \ref{conj2} above) in $\mathbb{T}^2$.

Notice that all these results are neither probabilistic nor semiclassical, which is in contrast with part of the recent literature.  At the same time, these observations together with the general conjectures stated above suggest their truth in higher dimensions as well. However, this is not the case (cf. Theorem \ref{cx}). 

We will next describe a distribution symmetry result in all dimensions, but only for a special class of trigonometric polynomials which does not include all eigenfunctions; this will be later complemented with counterexamples exhibiting its sharpness. To state it, we introduce the following class of trigonometric polynomials: 

\[f(x)=\sum_{i=1}^n(a_n\sin(2\pi \nu_i\cdot x)+b_n\cos(2\pi \nu_i\cdot x))\]
where $x\in\mathbb{T}^n$, $a_i,b_i\in\mathbb{R}$ and the $\nu_i\in\mathbb{Z}^n$ are linearly independent vectors. We shall denote this class of functions by $\mathcal{S}(\mathbb{T}^n)$. The main result we obtain for this class is 

\begin{theorem}[Distributional symmetry of functions in $\mathcal{S}$]\label{main2}
The distribution function $d\lambda$ of a non constant function $f\in\mathcal{S}(\mathbb{T}^n)$ is symmetric around $s=0$. In particular, it satisfies sign equidistribution and its global extrema coincide in absolute value.
\end{theorem}

Notice that this is neither contained nor implied by our previous results. It complements Theorem \ref{main} in dimension two, and analogous results regarding sign equidistribution, global maxima and $L^p$ norm symmetry follow for this class as well.

The linear independence hypothesis can not be dropped in any dimension as the example 
\[f(x)=\sin(x)+\cos(2x)\]
and small perturbations of it readily shows.

A more subtle counterexample within the class of eigenfunctions also exists, showing that, in general, Conjecture \ref{conj1} does not hold true. 

\begin{theorem}\label{cx}
The function $g(x,y,z)$ given by
\[\sin(x+y)-\cos(y-z)-\sin(x+z)\]
satisfies $-\Delta\psi=2\psi$ in $\mathbb{T}^3$ and it is negative for at least $52\%$ of the volume.
\end{theorem}

The proof of this is a computer assisted argument requiring about a billion computations.

The paper is organized as follows. In the next section we provide a proof of Theorem \ref{main} \color{blue}, \color{black} of which Theorems \ref{sign}, \ref{maxmin} and Corollary \ref{lp} are immediate consequences. The proof method applies to more general two dimensional tori. In Section \ref{higher} we show that it does not apply to higher dimensional tori.  In Section \ref{mainsection} we give three proofs of Theorem \ref{main2} on trigonometric polynomials in the class $\mathcal{S}$. They have different flavors although two of them are essentially equivalent. Finally, we devote Section \ref{sectioncx} to explain how to perform the computer assisted proof of Theorem \ref{cx}.

\section{Proof of theorem \ref{main}}

We will divide the proof in a number of cases depending on the eigenvalue. This is based on the following tricotomy: since an eigenvalue should be a sum of two squares, say,
\[\lambda=n^2+m^2\,,\]
then the following combinations of parities arise naturally: 
\begin{itemize}
\item[(a)] Both $n$ and $m$ are even iff $\lambda\equiv 0\mod 4$.
\item[(b)] The pair $n$ and $m$ have different parity iff $\lambda\equiv 1 \mod 4$.
\item[(c)] Both $n$ and $m$ are odd iff $\lambda\equiv 2\mod 4$.
\end{itemize}

Notice that the case $\lambda\equiv 3\mod 4$ is not possible in dimension two. The first two cases are easy to handle. 

For instance, (a) reduces to the latter. Indeed, if $\psi(x)$ is an eigenfunction with eigenvalue $\lambda\equiv 0\mod 4$ it is easy to observe that $\psi(x/2)$ will be an eigenfunction with eigenvalue $\lambda/4$ and the proportion of area where it is positive or negative is preserved. We might apply this procedure until we hit the cases (b) or (c).  

In case (b), we claim that the eigenfunction satisfies the following identity
\[\psi(x,y)=-\psi\left(x+\frac{1}{2},y+\frac{1}{2}\right)\]
which proves that the set where it is positive is a translation of the set where it is negative and viceversa. Since translations preserve volume the result follows. To see the claim, observe that since we can express
\begin{equation}\label{*}
\psi(x,y)=\sum_{\nu} a_{\nu}\exp(2\pi i \nu_1x+2\pi i\nu_2y)
\end{equation}
where the sum extends over the set of $\nu\in\mathbb{Z}^2$ such that $\nu_1^2+\nu_2^2=\lambda$. The proof in this case ends observing that the sum $\nu_1+\nu_2$, being odd, immediately implies the claimed functional identity.

Finally, in case (c) we know both coordinates are odd so it is enough to translate by $\frac{1}{2}$ in one of them to get the same functional identity. This concludes the proof as the translational symmetry is an involution and isometry, therefore: it interchanges the sets with different sign, global extrema, and the level sets in general. 

\section{Extensions of the argument and obstructions}\label{higher}

The same method of proof applies to other situations as well. For instance,  Theorem \ref{main} might be generalized to the following
 
\begin{theorem}[Symmetry]
Let $\Lambda \subset \mathbb{R}^{2}$ be a lattice and construct the torus $T=\mathbb{R}^2/(2 \pi \Lambda)$. Given a non constant real eigenfunction $\psi$ on $T$, the distribution function $d\lambda$ is symmetric around $s=0$.
\end{theorem}
The proof is a straightforward generalization of the one presented above, taking into account that the eigenfunctions are combinations of $\exp(i k \cdot x)$ with $k=n v_1+ m v_2$, where $(v_1, v_2)$ are the generators of the dual lattice $\Lambda^{*}$. Indeed, case (b) is dealt with using the translation $(x,y)\mapsto(x,y)+u$, with $u$ being the unique vector such that $v_1 \cdot u= v_2 \cdot u= \frac{1}{2}$, and the rest of the argument follows mutatis mutandis.
We leave further details to the reader. 

Likewise, the method of proof immediately provides 

\begin{theorem}\label{oddfrequencies}
Distributional symmetry holds for real eigenfunctions in $\mathbb{T}^3$ with eigenvalue $\lambda\not\equiv 2\mod 4$ or functions in any $\mathbb{T}^n$ supported in odd frequencies in any torus $\mathbb{T}^n$.
\end{theorem}

Notice that this is more general as it does not restrict to eigenfunctions. One might be led to believe that the exceptions to the theorem could be handled using a different argument. This is false, as the counterexample constructed in Theorem \ref{cx} shows. As the proof of this will be computer assisted we believe it is suitable to provide the reader beforehand with an elementary argument showing that there are no semiintegral translation $T:x\mapsto x+\frac{1}{2}v$ with $v\in\mathbb{Z}^3$ such that the identity $\psi(x)=-\psi(Tx)$ holds, obstructing the extension of the arguments above. In fact, this played a crucial role in finding the counterexample itself.

Indeed, suppose that the Fourier transform of $\psi$ is  supported at least in the frequency vectors $\nu_1^*=(1,1,0)$, $\nu_2^*=(1,-1,0)$, $\nu_3^*=(0,1,1)$, $\nu_4^*=(0,1,-1)$, $\nu_5^*=(1,0,1)$ and $\nu_6^*=(1,0,-1)$. It is evident that one needs the scalar product of any of these  with  $v=(x,y,z)$ to be an odd number which forces at least one of the entries to be odd. Due to the symmetry it is enough to suppose $x$ is odd. Let us now suppose  that $y$ is odd as well. The four possible scalar products of such a $v$ against the vectors $(1,1,0)$ and $(1,-1,0)$ show a contradiction, so $y$ must be even.  The same argument would apply to $z$. But this is a contradiction as the vector product with  $(0,1,1)$ would be even.

\section{Proof of Theorem \ref{main2}}\label{mainsection}

We present three proofs: one analytic, the other geometric and the last one purely algebraic. We find it appropriate to present the analytic proof, although it is not the most elementary, for a couple of reasons. First because it might be of interest in itself, and we hope the reader might find applications of the argument to different problems. Secondly, it was our first proof and the way we discovered the result originally.

\subsection{Analytic proof:}\label{anproof}

\color{black}
 
The proof will be based on the following

\begin{lemma}\label{MS}
Any continuous odd function can be uniformly approximated by  linear combinations of $x^{2k+1}$ for $k=0,1,\ldots$ in any symmetric interval $[-L,L]$.
\end{lemma}

The proof can be found in the appendix and is based on a simple application of the Muntz-Sz\'asz theorem (cf. \cite{FN} p. 114). We are now ready to fix $f\in\mathcal{S}$. Let us define $L=\|f\|_{\infty}$. Let us approximate the function $\textrm{sign}(x)$ by a continuous odd function $\sigma(x)$ such that $|\sigma(x)|\leq 1$ and $\sigma(x)\neq\textrm{sign}(x)$ for $|x|\leq \eta$ while $\sigma(x)=\textrm{sign}(x)$ for $|x|>\eta$. This implies
\[\int_{\mathbb{T}^n}\textrm{sign}(f(x))dx=\int_{\mathbb{T}^n}\sigma(f(x))dx+O(\operatorname{vol}(\{x:|f(x)|<\eta\})).\]
The latter is a thin neighbourhood of the nodal set $Z(f)$ of $f$. In fact, by continuity of $f$ and compactness it follows that for any $\delta>0$ there exist $\eta=\eta(\delta,f,\lambda)>0$ such that
\[\{f(x)\leq\eta\}\subseteq Z(f)+B_{\delta}.\]
The size of which can be shown to be $o(1)$ as $\eta$ tends to zero. 

Let us now estimate the integral in the right hand side. Given any $\epsilon>0$ an application of the Lemma \ref{MS} shows that
\[\int_{\mathbb{T}^n}\sigma(f(x))dx=\int_{\mathbb{T}^n}\sum_{k=0}^Na_k f(x)^{2k+1}dx+O(\epsilon).\]
The proof will end if we can proof that each integral
\[\int_{\mathbb{T}^n}f(x)^{2k+1}dx\]
vanishes. Indeed, in such a case one might let $\epsilon, \eta\rightarrow 0$ and putting together both estimates we would be done. To complete the proof then let us recall that $f$ can be expressed as in equation \ref{*}. Unfolding  the $(2k+1)$-fold product, the only terms that survive are those for which the frequencies satisfy
\begin{equation}\label{identity}
0=\sum_{j=0}^{2k+1}\nu_{i(j)}=\sum_{i=1}^nA_i\nu_i
\end{equation}
which can not happen because the $\nu_i$ are linearly independent. This concludes the proof of the sign equidistribution analogue for trigonometric polynomials in the class $\mathcal{S}$. 

Let us point out before continuing that the proof breaks down in general as the linear combination 
\[\nu^*_1-\nu_4^*-\nu_5^*=(1,1,0)-(0,1,-1)-(1,0,1)=0\]
shows (cf. Section \ref{higher}, last paragraph, for the definition of $\nu^*$).

The same method of proof provides the generalization

\begin{proposition}
For any $k\in\mathbb{N}$ the identity
\[\int_{\mathbb{T}^2}f(x)^{2k}\operatorname{sign}(f(x))dx=0\]
holds.
\end{proposition}

We have proved $k=0$ above. We leave details to the reader. This implies the following identity
\[\int_0^Ls^{2k}d\operatorname{vol}(\{f(x)>s\})=\int_0^Ls^{2k}d\operatorname{vol}(\{f(x)<-s\}).\]
This and a straightforward variation of the uniqueness results known for the moment problem\footnote{One only needs to replace the application of Weierstrass' theorem in Chapter II Section 6 from \cite{W} by a version of the M\"untz-Sz\'asz theorem.} imply that the distribution functions
\[\{f(x)>s\}\textrm{ and }\{-f(x)>s\}\]
 coincide for $s>0$ which concludes the proof. 
 
\subsection{Geometric proof}
Without loss of generality let $\{\nu_i\}_{i=1}^n$ be a basis of $\mathbb{R}^n$. It is possible to change to the canonical basis by a linear map whose entries are rational. This induces   a transformation on the torus which distorts the volume uniformly (the same distortion at every point). As a consequence the distribution functions of our function and the function in the new coordinates are proportional. The latter satisfies the distribution symmetry conjecture since there is a translation $T$ such that $g(Tx)=-g(x)$. 

\subsection{Algebraic proof}
Given $f$ expressed as in equation \ref{*}, let $A$ be the $n \times n$ matrix whose $i$th column is given by the frequency vector $\nu_i$; equivalently, let $A_{ji}:=(\nu_{i})_j$. Observe that, denoting by $\{e_i\}_{i=1}^{n}$ the standard basis in $\Rn$, we have $A e_i=\nu_i$.

Since the frequencies $\{\nu_i\}_{i=1}^{n}$ are linearly independent, the matrix $A$ has an inverse $A^{-1}$. In general, the determinant of $A$ will be (in absolute value) bigger than $1$, and thus the entries of $A^{-1}$ will not be integers, but rational numbers (cf. end of Section 4 above).

Denote by $(A^{-1})^{t}$ its transpose. For a sufficiently large $N$, the matrix $B:=N (A^{-1})^{t}$ has integer entries and defines a map $\Phi: \mathbb{T}^{n} \rightarrow \mathbb{T}^{n}$. It is easy to see that this map is a covering with degree equal to the determinant of $B$; locally, it is a diffeomorphism with constant Jacobian, i.e, with uniform volume distortion. Thus the distribution function of $f$ is proportional to that of the function
\[
g(x):=f(B x)=\sum_{i=1}^{n}  a_i \cos(B^{t}\nu_i \cdot x)+b_i \sin (B^{t}\nu_i \cdot x).
\]
On the other hand, $B^{t}\nu_i=N e_i$, so
\[
g(x)=\sum_{i=1}^{n}  a_i \cos (2 \pi N x_i)+b_i \sin (2 \pi N x_i)\,.
\]
Therefore $g$ has a translational antisymmetry $g(x+w)=-g(x)$, for the vector $w=\frac{1}{2N}(1, 1,...,1)^{t}$. Since translations preserve volume, we have distribution symmetry for $g$, and hence also for $f$.

A more succinct variant of the previous one, only without a geometric interpretation. With the previous notations, define
\[
u:=\frac{1}{2}\sum_{i=1}^{n} (A^{-1})^{t} e_{i} \,.
\]

This vector verifies, for any $i=1,...,n$, the identity $\nu_i \cdot u=\frac{1}{2}$. Thus we conclude that the function $f$ is antisymmetric under translations by $u$, i.e

\[
f(x+u)= \sum_{i=1}^{n}  a_i \cos(2 \pi \nu_i x+\pi)+b_i \sin (2 \pi \nu_i x+\pi)=-f(x) \,,
\]

and the result follows.

\section{Proof of Theorem \ref{cx}}\label{sectioncx}

We can estimate the volume of $(x,y,z)\in[0,1]^3$ such that 
\[g(x,y,z)=\sin(2\pi(x+y))-\cos(2\pi(y-z))-\sin(2\pi(x+z))<0\]
from below numerically. Indeed, it is easy to check that the gradient of $g$ is uniformly bounded by $6\pi$, this shows that  the error $g(x)-g(y)$, is bounded by $e=6\sqrt{3}\pi L$ for any point $y$ within a cube of sidelength $L$ centered at $x$. This allows us to  reduce the problem to a counting argument as it will be enough to  find the proportion of cubes where $g$ is negative. As a consequence, if the mesh length $L$ of a grid $G$ is small enough, the volume where $g$ is negative is bounded from below by the proportion of the centers of cubes $x\in G$ such that
\[g(x)<-e.\]
The above argument provides a criteria to bound this proportion from below that can be done computer assistedly for a large grid. However, the machine will commit rounding errors when computing the elementary function $g$ at any $x\in G$. To address this issue rigorously we will check the condition
\[\bar{g}(x)<-1.1e,\]
where we denote by $\bar{g}(x)$ the machine output of the calculation of $g(x)$.   Indeed, the computation to check the above inequality requires to compute $g$ at every $x\in G$ with a certain accuracy, that, say, should be smaller than $0.05e$. As a consequence, for any $y$ in the cube centered at $x$
\[g(y)=g(y)-g(x)+g(x)-\bar{g}(x)+\bar{g}(x)<e+0.05e-1.1e<0.\]
Of course, we are forced to use an approximation of $\pi$. Let the approximation be $\bar{\pi}$. The mean value theorem implies that the error commited when we compute each trigonometric function (at any fixed $x\in[0,1]^3$) is bounded by $4\cdot |\pi-\bar{\pi}|$. (Let us remark that the mean value theorem is applied with respect to the variable $\bar{\pi}$.) In our computations we can use an approximation $\bar{\pi}$ with accuracy $ 10^{-20}$. This adds an error corresponding to the computation of the three trigonometric functions with $\bar{\pi}$ which together is then bounded by $12\cdot 10^{-20}$. The error commited computing 
\[g'(x,y,z)=\sin(2\bar{\pi}(x+y))-\cos(2\bar{\pi}(y-z))-\sin(2\bar{\pi}(x+z))<0\]
 in $C$++ is bounded by, at most, $3\cdot 10^{-6}$. Finally, the computer will need to add this three numbers, which adds rounding errors $10^{-5}$  twice. The computation will be reliable if $11\cdot 10^{-6}<0.05e$ holds. 

Let us emphasize that this is a crude estimate on the error and using more precission arithmetic one can do much better but as we will see below for this problem we will not need to go deeper as all we are approximating is the sum of three trigonometric functions.

A simple algorithm storing in a variable $M$ the number of points $x\in G$ that satisfy the condition 
\[g(x)<-1.1e\]
allows to compute the proportion $M/L^3$. Due to the unusual number of computation required this is achieved through a computer assisted approach. We supplemented this argument with an algorithm implemented in $C$++ which in the particular case of a grid with $N^3=(2^{7})^3$ points gives
\[\frac{M}{N^3}=\frac{1123200}{2097152}\]
which is roughly an estimated $53.5\%$ of points $x\in G$  satisfying $g(x)<-1.1e<0$. In this case the mesh sidelength is $L=1/N=2^{-7}$ which implies the maximal error $e$ committed within each box centered at those grid points is bounded above by $0.26$. As a consequence the proportion of points provides a lower bound estimate for the points where $g(x)<0$ proving the result. Similarly one can conclude that the size of the area positiveness is at least $34.3\%$. This is far from giving a close estimate. Using about a billion nodes instead ($(2^{10})^3$ nodes) one can get the lower estimates $59.3\%$ for the set of points where $g$ is negative and $39.1\%$ for its complementary. This time the uncertainty of sign reduces to less than $2\%$ of the total volume.

The  precision of any machine for the required calculation is within $10^{-6}$, an accuracy that makes the above argument reliable and allows us to avoid a more rigorous study of the cumulative errors the algorithm might have committed.

\section{Appendix}

\textsc{Proof of Lemma \ref{MS}:} the Muntz-Sz\'asz theorem implies that any continuous function $f$ in $[0,1]$ can be uniformly approximated by linear combinations of $1$ and $x^{2k+1}$ with $k=0,1,\ldots$. We might use a rescaling to let $L=1$ without loss of generality. The oddness in our case implies that $f(0)=0$ from which it follows that for any $\epsilon$ there exist  coefficients $a_n$ such that
\[f(x)=a_0+\sum_{k=0}^N a_k x^k+O(\epsilon)\]
  and   evaluating at $x=0$ it follows that $a_0=O(\epsilon)$ too and we might forget about it for free. The rest follows by odd symmetry.

\section{\textbf{Acknowledgments}}

The authors are grateful to the Fields Institute Hydrodynamics Program and the University of Toronto for their support. F.T.L gratefully acknowledges support from the Max-Planck Institute for Mathematics. The authors would like to thank the referees for their suggestions.

\end{document}